\input amstex
\documentstyle{amsppt}
\magnification=\magstep1
 \hsize 13cm \vsize 18.35cm \pageno=1
\loadbold \loadmsam
    \loadmsbm
    \UseAMSsymbols
\topmatter
\NoRunningHeads
\title A note on   $q$-Euler numbers and polynomials
\endtitle
\author
  Taekyun Kim
\endauthor
\abstract
 The purpose of this paper is to construct $q$-Euler numbers and polynomials by using $p$-adic
 $q$-integral equations on $\Bbb Z_p$. Finally, we will give some
 interesting formula related to these $q$-Euler numbers and
 polynomials.
\endabstract
\endtopmatter

\document

{\bf\centerline {\S 1. Introduction}}

 \vskip 20pt

The usual Bernoulli numbers are defined by
$$
\sum_{k=0}^\infty B_k \dfrac{t^k}{k!} =\dfrac{t}{e^t -1},
$$
which can be written symbolically as $e^{Bt} =\frac{t}{e^t -1}$,
interpreted to mean $B^k$ must be replaced by $B_k$ when we expand
on the left.  This relation can also be written $e^{(B+1)t} -
e^{Bt} =t$, or, if we equate power of  $t$,
$$
B_0 =1, \quad (B +1)^k -B_k = \cases
1 & \text{if\   $k=1$}\\
0 & \text{if \  $k>1$,}
\endcases
$$
where again we must first expand and then replace $B^i$ by $B_i$,
cf. [6,8,9,10,11].

Carlitz's $q$-Bernoulli numbers $\beta_k$ can be determined
inductively by
$$\eqalignno{ &
\beta_0 =1,\quad  q(q\beta +1)^k -\beta_k = \cases
1 & \text{if\   $k=1$}\\
0 & \text{if \  $k>1$,}
\endcases &(1)}
$$
with the usual convention about replacing $\beta^i$ by
$\beta_{i,q}$(see [1,2,3,4,5,12]).

Carlitz also defined $q$-Euler numbers and polynomials as
$$\eqalignno{ &
H_0 (u;q)=1, \quad (qH +1)^k -u H_k (u; q) =0 \quad \text{ for }
k\geq 1, &(2)}$$
  where $u$ is a complex number with $|u|>1$; and for $k\geq 0$,
  $$H_k (u, x ; q) = (q^x H +[x]_q )^k , \quad \text{ cf. }
  [2,3,4,9,10,11,12],
  $$
where $[x]_q$ is defined by $[x]_q = \frac{1-q^x}{ 1-q}$.

It was known that the ordinary Euler polynomials are defined by
$$\eqalignno{ &
\dfrac{2}{e^t +1}e^{xt} = \sum_{n=0}^\infty E_n (x)
\dfrac{t^n}{n!}, \quad (|t|<\pi ),&(3)}
$$
we note that $E_n =E_n (0)$ are called $n$-th Euler numbers, cf.
[6,7,8].

Let $p$ be a fixed odd prime, and let $\Bbb C_p$ denote the
$p$-adic completion of the algebraic closure of $\Bbb C_p$. For
$d$ a fixed positive integer $(p,d)=1$, let
$$\split
& X=X_d = \lim_{\overleftarrow{N} } \Bbb Z/ dp^N \Bbb Z ,\cr & \
X_1 = \Bbb Z_p , \cr  & X^\ast = \underset {{0<a<d p}\atop
{(a,p)=1}}\to {\cup} (a+ dp \Bbb Z_p ), \cr & a+d p^N \Bbb Z_p =\{
x\in X | x \equiv a \pmod{dp^n}\},\endsplit$$ where $a\in \Bbb Z$
lies in $0\leq a < d p^N$.

The $p$-adic absolute value in $\Bbb C_p$ is normalized so that
$|p|_p = \frac{1}{p}$. Let $q$ be variously considered as an
indeterminate a complex number $q\in\Bbb C$, or a $p$-adic number
$q\in\Bbb C_p$. If $q\in\Bbb C$, we always assume $|q|<1$. If
$q\in\Bbb C_p$, we always assume $|q-1|_p <p^{-1/p-1}$, so that
$q^x =\exp (x\log q)$ for $|x|_p \leq 1$.

Throughout this paper we use the notation:
$$
[x]_q  = \dfrac{q^x -1}{q -1} = 1+ q + q^2 +\cdots + q^{x-1}.
$$

We say that $f$ is uniformly differentiable function at a point $a
\in\Bbb Z_p $ and denote this property by $f\in UD(\Bbb Z_p )$, if
the difference quotients
$$
F_f (x,y) = \dfrac{f(x) -f(y)}{x-y}
$$
have a limit $l=f^\prime (a)$ as $(x,y) \to (a,a)$.

For $f\in UD(\Bbb Z_p )$, let us start with expression

$$\eqalignno{ & \dfrac{1}{[p^N ]_q} \sum_{0\leq j < p^N} q^j f(j) =\sum_{0\leq j < p^N} f(j)
\mu_q (j +p^N \Bbb Z_p ),&(4) }
$$
representing $q$-analogue of Riemann sums for $f$, cf. [5]. The
integral of $f$ on $\Bbb Z_p$ will be defined as limit ($n \to
\infty$) of those sums, when it exists. The $q$-Volkenborn
integral of function $f\in UD(\Bbb Z_p )$ is defined by
$$ \eqalignno{ &
I_q (f) =\int_{\Bbb Z_p }f(x) d\mu_q (x) = \lim_{N\to \infty}
\dfrac{1}{[p^N ]_q} \sum_{0\leq x < p^N} f(x) q^x ,\quad
\text{(see [5] )}. &(5)}
$$
From the definition of $[x]_q$, we derive
$$
[x]_{-q} = \dfrac{1-(-q)^x}{ 1+q} = 1-q +q^2 -q^3 + \cdots
+(-1)^{x-1} q^{x-1} .
$$
In [4,5], it was known that

$$\eqalignno{ &\int_{\Bbb Z_p} [x+y]_q^m d\mu_q (y)
=\beta_{m,q} (x),&(6) }
$$
where $\beta_{m,q} (x)$  are called Carlitz's $q$-Bernoulli
polynomials.  By (5), we easily see that

$$\eqalignno{ & I_{-q} (f) = \int_{\Bbb Z_p} f(x) d\mu_{-q} (x)
= \dfrac{[2]_q}{2} \lim_{N\to \infty }\sum_{x=0}^{p^N -1} f(x)
(-1)^x q^x . &(7)}$$

The purpose of this note is to construct $q$-Euler numbers which
can be uniquely determined by
$$\eqalignno{ & q(qE_q +1 )^n + E_{n,q} =\cases
[2]_q & \text{if\   $n=0$}\\
0 & \text{if \  $n>0$,}
\endcases &(8)}
$$
with the usual convention about replacing $E_q^i$ by $E_{i,q}$.
From these numbers, we will derive some interesting formulae.

\vskip 20pt

{\bf\centerline {\S 2. A note on $q$-Euler numbers and
polynomials}}

 \vskip 20pt

From (7), we derive formulae as follows:
$$\eqalignno{ &
qI_{-q} (f_1 ) + I_{-q} (f) =[2]_q f(0), &(9)}
$$
where $f_1 (x)$ is translation with $f_1 (x) =f(x+1)$.

 Let $f(x) =
e^{[x]_q t}$. Then we see that
$$\eqalignno{ & q \int_{\Bbb Z_p} e^{[x+1]_q t} d\mu_{-q} (x) +
\int_{\Bbb Z_p} e^{[x]_q t} d\mu_{-q} (x) = [2]_q  . &(10)}
$$

First, we consider the following integral :

$$\eqalignno{
\int_{\Bbb Z_p} e^{[x]_q t} d\mu_{-q} (x) &= \sum_{n=0}^\infty
\int_{\Bbb Z_p} [x]_q^n d\mu_{-q} (x) \dfrac{t^n}{n!}\cr &=
\sum_{n=0}^\infty \left( \dfrac{1}{1-q} \right)^n \sum_{l=0}^n
\binom{n}{l} (-1)^l \int_{\Bbb Z_p} q^{xl} d\mu_{-q} (x)
\dfrac{t^n}{n!}\cr &= [2]_q \sum_{n=0}^\infty \left(
\dfrac{1}{1-q} \right)^n \sum_{l=0}^n \binom{n}{l}  (-1)^l
\sum_{m=0}^\infty   (-1)^m q^{(l+1)m} \dfrac{t^n}{n!}\cr &= [2]_q
\sum_{m=0}^\infty (-1)^m q^{m} e^{[m]_q t}= \sum_{n=0}^\infty
E_{n,q} \dfrac{t^n}{n!}. &(11)}
$$
By (11), we easily see that
$$\eqalignno{ &
\int_{\Bbb Z_p} [x]_q^n d \mu_{-q} (x) = E_{n,q} .  &(12)}
$$

 Note that
 $$
\lim_{q \to 1} E_{n,q} (x) = E_n =\int_{\Bbb Z_p} x^n d \mu_{-q}
(x) , \quad \text{ see}[7].$$

By the same method, we note that

$$\split
\int_{\Bbb Z_p} e^{[x+y]_q t} d\mu_{-q} (y) &= \sum_{n=0}^\infty
\int_{\Bbb Z_p}[x+y]_q^n d \mu_{-q} (y) \dfrac{t^n}{n!}\cr
 &= \sum_{n=0}^\infty \left(\dfrac{1}{1-q} \right)^n \left(
 \sum_{l=0}^n \binom{n}{l} (-1)^l q^{lx} \int_{\Bbb Z_p}
 q^{yl}d\mu_{-q} (y)
\right) \dfrac{t^n}{n!}\cr
 &= [2]_q \sum_{n=0}^\infty
\left(\dfrac{1}{1-q} \right)^n \left(
 \sum_{l=0}^n \binom{n}{l} \dfrac{(-1)^l q^{lx}}{1+q^{l+1}}\right) \dfrac{t^n}{n!}\cr
&= [2]_q \sum_{m=0}^\infty (-1)^m q^m
 \sum_{n=0}^\infty [m+x]_q^n  \dfrac{t^n}{n!}\cr
&= [2]_q \sum_{m=0}^\infty (-1)^m q^m
 e^{[x+m]_q t}\cr
 &=\sum_{m=0}^\infty E_{n,q} (x) \dfrac{t^n}{n!}.
\endsplit
$$

Thus, we have
$$\eqalignno{ &\int_{\Bbb Z_p} [x+y]_q^n  d \mu_{-q} (y ) =
 E_{n,q} (x) ,\quad  n \geq 0. &(13)}
 $$
Note that
$$
\lim_{q \to 1} E_{n,q} (x) = E_n (x) = \int_{\Bbb Z_p} (x+y)^n d
\mu_{-1} (y).
$$

Since
$$\eqalignno{ & [x+1]_q = \dfrac{1-q^{x+1}}{1-q} =\dfrac{1-q}{1-q}
+ \dfrac{1-q^x}{1-q}q = 1 +q [x]_q . &(14)}
$$
From (10), (12) and (14), we derive
$$\eqalignno{ [2]_q & =\sum_{n=0}^\infty  \left(
q \int_{\Bbb Z_p} [x+1]_q^n d \mu_{-q} (x) +  \int_{\Bbb Z_p}
[x]_q^n d \mu_{-q} (x)
 \right) \dfrac{t^n}{n!}\cr
 &=\sum_{n=0}^\infty  \left(
q \sum_{l=0}^n \binom{n}{l} q^l \int_{\Bbb Z_p} [x]_q^l d \mu_{-q}
(x) + \int_{\Bbb Z_p} [x]_q^n d \mu_{-q} (x)
 \right) \dfrac{t^n}{n!}\cr
 &=\sum_{n=0}^\infty  \left(
q \sum_{l=0}^n \binom{n}{l} q^l E_{l,q} + E_{n,q}
 \right) \dfrac{t^n}{n!}\cr
 &=\sum_{n=0}^\infty  \left(
q (q E_q +1 )^n + E_{n,q}
 \right) \dfrac{t^n}{n!}, &(15)}
$$
with the usual convention about replacing $E_q^i$ by $E_{i,q}$.

By (15), we easily see that
$$\eqalignno{ & q(qE_q +1)^n +E_{n,q} =
\cases
[2]_q & \text{if\   $n=0$}\\
0 & \text{if \  $n>0$,}
\endcases &(16)}
$$
with the usual convention about replacing $E_q^i$ by $E_{i,q}$.
When we compare Eq.(16) and Eq.(1), the Eq.(16) seems to be
interesting formula. In particular, these numbers seem to be new,
which are different than Carlitz's $q$-Euler numbers. From (7), we
derive
$$\eqalignno{&
q^n I_{-q} (f_n ) + (-1)^{n-1} I_{-q} (f) = [2]_q\sum_{l=0}^{n-1}
(-1)^{n-1-l} q^l f(l), &(16-1)}
$$
where $n\in\Bbb N$, $f_n (x) =f(x+n)$. When $n$ is an odd positive
integer, we note that
$$\eqalignno{&
q^n I_{-q} (f_n ) + I_{-q} (f) =[2]_q \sum_{l=0}^{n-1} (-1)^l q^l
f(l). &(17)}
$$
From (17), we derive
$$\eqalignno{&
q^n \int_{\Bbb Z_p} [x+n]_q^m d\mu_{-q} (x) + \int_{\Bbb Z_p}
[x]_q^m d\mu_{-q} (x) =[2]_q \sum_{l=0}^{n-1} (-1)^l q^l [l]_q^m ,
&(18)}
$$
where $n$ is an odd positive integer. By (12) and (13), we easily
see that
$$\eqalignno{& [2]_q \sum_{l=0}^{n-1} (-1)^l q^l [l]_q^m = q^n
E_{m,q} (n) + E_{m,q}, &(19)}
$$
where $n$ is an odd positive integer.  If $n$
 is an even integer, then we have in Eq.(16-1) that

$$\eqalignno{& q^n \int_{\Bbb Z_p} [x+n]^m_q d \mu_{-q} (x) - \int_{\Bbb Z_p } [x]_q^m d \mu_{-q} (x)
=[2]_q \sum_{l=0}^{n-1} (-1)^{l-1} q^l [l]_q^m . &(20)}
$$

From (12), (13) and (20), we derive
$$\eqalignno{& q^n E_{m,q}(n) -E_{m,q}=
[2]_q  \sum_{l=0}^{n-1} (-1)^{l-1} q^l [l]_q^m , &(21)}
$$
where $n$ is a positive even integer. It seems to be interesting
to compare (19) and (21).

\vskip 20pt

\quad Taekyun Kim

\quad EECS, Kyungpook National University, Taegu 702-701, S. Korea

\quad e-mail:\text{ tkim$\@$knu.ac.kr; tkim64$\@$hanmail.net}

\enddocument